\documentclass[12pt,a4paper]{article}
\usepackage{amsmath,amssymb,amsfonts}
\def\Bbb{\mathbb}

\title{\bf The largest values of Dedekind sums}

\author{Kurt Girstmair}
\date{}

\makeatletter
\let\@@maketitle=\maketitle
\def\maketitle{\def\thispagestyle##1{\relax}\@@maketitle}
\makeatother
\newtheorem{theorem}{Theorem}

\def\BE{\begin{equation}}
\def\EE{\end{equation}}
\def\BD{\begin{displaymath}}
\def\ED{\end{displaymath}}
\def\BA{\begin{array}}
\def\EA{\end{array}}
\def\BEA{\begin{eqnarray*}}
\def\EEA{\end{eqnarray*}}
\def\BI{\bibitem}

\def\Z{\Bbb Z}

\def\R{\Bbb R}

\def\phi{\varphi}
\def\EPS{\varepsilon}

\def\MB{\mbox}
\def\LD{\ldots}

\def\BQ{``}
\def\EQ{'' }

\def\MN{\medskip\noindent}
\def\STOP{\hfill$\Box$}

\def\DED{Dedekind }

\begin{document}
\maketitle

\begin{abstract}

\noindent
Let $s(m,n)$ denote the classical \DED sum, where $n$ is a positive integer and $m\in\{0,1,\LD, n-1\}$,
$(m,n)=1$. For a given positive integer $k$, we describe a set of at most $k^2$ numbers $m$ for which $s(m,n)$ may be
$\ge s(k,n)$, provided that $n$ is sufficiently large. For the numbers $m$ not in this set, $s(m,n)<s(k,n)$.

\end{abstract}

\section*{1. Introduction and results}

Let $m$ and $n$ be integers, $n\ne 0$ and $(m,n)=1$. The classical \DED sum $s(m,n)$ is defined by
\BD
   s(m,n)=\sum_{k=1}^{|n|} ((k/n))((mk/n))
\ED
where $((\LD))$ is the \BQ sawtooth function\EQ defined by
\BD
  ((t))=\left\{\begin{array}{ll}
                 t-\lfloor t\rfloor-1/2 & \MB{ if } t\in\R\smallsetminus \Z; \\
                 0 & \MB{ if } t\in \Z
               \end{array}\right.
\ED
(see, for instance, \cite[p. 1]{RaGr}).

In the present setting it is more
convenient to work with
\BD
 S(m,n)=12s(m,n)
\ED instead.
Since $S(m,-n)=S(m,n)$ and $S(m+n,n)=S(m,n)$, we obtain all \DED sums if we restrict $n$ to positive integers and $m$ to the range $0\le m< n$.
The general case, however, will be needed below (see (\ref{2.2})).

The original context of \DED sums is the theory of modular forms (see \cite{Ap2}). But these sums have also interesting applications in
connection with class numbers, lattice point problems, topology, and algebraic geometry (see  \cite{{At}, {Me}, {RaGr}, {Ur}}).
Starting with Rademacher \cite{Ra}, several authors have studied the distribution of \DED sums (for instance, \cite{{Br}, {Gi2}, {Hi}, {Va}}).
Whereas the arithmetic mean of the absolute values $|S(m,n)|$, $0\le m<n$, $(m,n)=1$, has order of magnitude $\log^2 n$ for $n$ tending to infinity (see \cite{GiSc}), large \DED sums $S(m,n)$ have order of magnitude $n$.

In this paper we study the largest values of \DED sums $S(m,n)$ for a given sufficiently large number $n$.

In 1956, Rademacher showed
\BE
\label{1.1}
S(1,n)>S(m,n) \MB{ for all }m, 2\le m\le n-1, (m,n)=1
\EE
(see \cite[Satz 2]{Ra}).
By the reciprocity law for \DED sums (see \cite[p. 5]{RaGr}),
\BE
\label{1.2}
  S(k,n)=-S(n,k)+ \frac nk+\frac kn+\frac 1{kn}-3,
\EE
we obtain
\BE
\label{1.3}
  S(1,n)=\frac{n^2-3n+2}{n}.
\EE
So this largest of
all \DED sums $S(m,n)$ equals
$n+O(1)$ for $n$ tending to infinity. Other large \DED sums are
$S(k,n)$
for a fixed integer $k>1$ and large numbers $n$. In fact, $S(k,n)=n/k + O(1)$, see (\ref{2.6}).  The main result of this paper is the following theorem.

\begin{theorem} 
\label{t1}
Let $k$ be a positive integer. For sufficiently large integers $n>k$ with $(k,n)=1$
and $m\in\{0,\LD,n-1\}$, $(m,n)=1$, we have
\BD
  S(m,n)\ge S(k,n)
\ED
only if $m$ has the form
\BE
\label{1.4}
  m=\frac{nc+q}d,
\EE
with
\BE
\label{1.5}
 d\in\{1,\LD,k\},
   c\in\{0,1,\LD,d-1\},(c,d)=1, q\in\{1,\LD, \lfloor k/d\rfloor\}.
\EE

\end{theorem} 

\MN
{\em Remarks.} 1. The proof of Theorem \ref{t1} shows
\BD
   S(m,n)=\frac n{dq}+O(1)
\ED
for each of the numbers (\ref{1.4}) in question. Since $dq\le k$, we see that $S(m,n)>S(k,n)$ whenever $dq<k$, whereas $S(m,k)\ge S(m,n)$ may hold if
$dq=k$. The proof of Theorem \ref{t1} also gives
\BD
  S(m,n)\le \frac n{k+1}+O(1)
\ED
for all numbers $m$ not of the form described by (\ref{1.4}) and (\ref{1.5}).

2. It is easy to see that there are at most
\BD
  \sum_{d=1}^k\phi(d)\frac kd.
\ED
numbers $m$ as described by (\ref{1.4}), (\ref{1.5}).
This bound is $\le k^2$ or, more precisely, $=6k^2/{\pi^2}+O(k\log k)$
(see \cite[p. 70, Ex. 5]{Ap}). In most cases, however, this is only a rough upper bound.

\MN
{\em Examples.} 1. For $k=3$, the numbers $m$ described by (\ref{1.4}), (\ref{1.5}) are $m=1,2,3$ (with $d=1, c=0$, $q=1,2,3$),
$m=(n+1)/2$ (with $d=2$, $c=1$, $q=1$), and $m=(n+1)/3, (2n+1)/3$ (with $d=3$, $c=1,2$, $q=1$). If $(n+1)/2$ is an integer, $n$ must be odd. Then $S((n+1)/2, n)=S(2,n)$, since $(n+1)/2$ is the inverse
of $2$ mod $n$. The last two cases occur only if $n\equiv 2$ mod 3 and $n\equiv 1$ mod 3, respectively. In each of these cases, $S(m,n)=S(3,n)$.

2. This example might suggest that for the numbers $m$ given by (\ref{1.4}), (\ref{1.5}) the \DED sums $S(m,n)$ take one of the values
$S(j,n)$, $j\in\{1,\LD,k\}$. This, however, is not true. Indeed, let $k=6$, $d=3$, $c=1$, $q=2$, and $n\ge 7$, so $m=(n+2)/3$. Since $m$
must be an integer, we require $n\equiv 1$ mod $3$. Because $k=6$, $n$ must be odd, and so $n\equiv 1$ mod 6. Under this condition,
we obtain from (\ref{2.2}) below
\BD
  S(m,n)=\frac{n^2-14n+13}{6n},
\ED
whereas the reciprocity law (\ref{1.2})
yields
\BD
  S(6,n)=\frac{n^2-38n+37}{6n}.
\ED
Accordingly, $S(m,n)=S(6,n)+O(1)$, but always $S(m,n)>S(6,n)$.

\MN
All terms $O(1)$ in this paper can be transformed into explicit bounds. In this way, one may obtain results of Rademacher type (see (\ref{1.1})) for any given $k$. As an example,
we settle the case $k=2$ here.

\begin{theorem} 
\label{t2}
Let $n\ge 3$ be odd (hence $S(2,n)$ is defined). Then for all $m\in\{3,\LD,n-1\}$, $(m,n)=1$, different from $(n+1)/2$,
\BD
  S(2,n)>S(m,n).
\ED

\end{theorem} 

\section*{2. Proofs}

{\em Proof of Theorem \ref{t1}.} Put $l=2k+2$ and let $n>l$. We call a number $m\in\{1,\LD,n-1\}$, $(m,n)=1$, {\em ordinary} if, and only if,
for all $d\in\{1,\LD,l\}$ and all $c\in\{0,1,\LD,d\}$ with $(c,d)=1$,
\BD
 \left|\frac mn-\frac cd\right|\ge \frac l{nd},
\ED
i.e., each possible $q=md-nc$ satisfies $|q|\ge l$. Let $m$ be an ordinary number. By a theorem about Farey approximation (see \cite[p. 127, Th. 10.5]{Hu}),
there is a number $d\in\{1,\LD,l\}$ and a $c\in\{0,1,\LD, d\}$, $(c,d)=1$, such that
\BD
 \left|\frac mn-\frac cd\right|\le \frac 1{ld}.
\ED
If we choose $d$ and $c$ in this way, we have $|q|\le n/l$ for the above $q$. Altogether,
\BD
  l\le|q|\le \frac nl.
\ED
By \cite[Lemma 1]{Gi},
\BE
\label{2.2}
 S(m,n)=S(c,d)+\EPS S(r,q)+\frac n{dq}+\frac d{nq}+\frac q{nd}-3\EPS
\EE
where $r$ is some integer prime to $q$ and $\EPS\in\{\pm 1\}$ is the sign of $q$ (observe $q\ne 0$ since $n>d$).
Combined with (\ref{1.3}), this gives
\BD
 |S(m,n)|\le d+|q|+\frac n{d|q|}+\frac d{n|q|}+\frac{|q|}{nd}+3.
\ED
We observe $d\le l$ and $|q|\le n/l$.
Further, since $|q|\ge l$, we have $n/(d|q|)\le n/l$. The condition $|q|\ge l\ge d$ implies $d/(n|q|)\le 1/n$. From
$|q|\le n/l$ we obtain $|q|/(nd)\le 1/l$. Altogether,
\BE
\label{2.4}
|S(m,n)|\le l+\frac nl+\frac nl+\frac 1n+\frac 1l+3=\frac{2n}l + l+\frac 1n + \frac 1l+ 3=\frac n{k+1}+O(1).
\EE

Next we show
\BE
\label{2.6}
S(k,n)=\frac nk+O(1).
\EE
To this end we observe $S(n,k)\le S(1,k)$ and $S(n,k)\ge -S(1,k)$, by (\ref{1.1}). Then
the reciprocity law (\ref{1.2}), combined with (\ref{1.3}), gives
\BD
S(k,n)\ge \frac{n^2-(k^2+2)n+k^2+1}{kn} \MB{ and } S(k,n)\le \frac{n^2+(k^2-6k+2)n+k^2+1}{kn}.
\ED
This implies (\ref{2.6}). Moreover, (\ref{2.4}) and (\ref{2.6}) show
\BD
S(m,n)<S(k,n)
\ED
for large numbers $n$ and ordinary numbers $m$.

Now suppose that $m$ is not an ordinary number. Therefore, there is a $d\in\{1,\LD,l\}$ and a $c\in\{0,1,\LD,d\}$,
$(c,d)=1$, such that $q=md-nc$ satisfies $|q|<l$. From (\ref{2.2}) we obtain
\BE
\label{2.8}
 S(m,n)=\frac n{dq}+ O(1),
\EE
where
\BE
\label{2.10}
 |O(1)|\le d+|q|+\frac d{n|q|}+\frac{|q|}{nd}+3\le 2l+\frac{2l}n+3\le 2l+5
\EE
because $n\ge l$. Accordingly, if $d|q|\ge k+1$, then $|S(m,n)|<S(k,n)$ for large numbers $n$. Thus, the only numbers $m$ to be considered
are those with $d|q|\le k$. They are, however, only of interest if $q>0$, since, otherwise, $S(m,n)<0$ according to (\ref{2.8}).
But these numbers are just those described by (\ref{1.4}), (\ref{1.5}).
\STOP

\MN
{\em Proof of Theorem \ref{t2}.} According to (\ref{2.4}), we have, for $k=2$ and ordinary numbers $m$,
\BD
 |S(m,n)|\le \frac n3+O(1),
\ED
with $|O(1)|\le 6+1/3+1/6+3=19/2$ since $n\ge 3$. Therefore, if $n/3+19/2<S(2,n)=(n^2-6n+5)/2n$, then $|S(m,n)|<S(2,n)$. This is the case for $n\ge 75$.

On the other hand, if $m$ is not an ordinary number, (\ref{2.8}) and (\ref{2.10}) give
\BD
  S(m,n)=\frac n{dq}+ O(1),
\ED
with $|O(1)|\le 2l+2l/n+3\le 12+12/3+3=19$. If $d|q|=1$ or $d|q|=2$, then $S(m,n)\ge S(2,n)$ only for $m=1,2,(n+1)/2$. For $|q|\ge 3$,
$S(m,n)<S(2,n)$ as soon as $n/3+19<S(2,n)$. This is the case for $n\ge 132$. Accordingly, Theorem \ref{t2} must be checked only for $n\le 131$, where it is also true.
\STOP


\vspace{0.5cm}
\noindent
Kurt Girstmair            \\
Institut f\"ur Mathematik \\
Universit\"at Innsbruck   \\
Technikerstr. 13/7        \\
A-6020 Innsbruck, Austria \\
Kurt.Girstmair@uibk.ac.at

\end{document}